\pdfoutput=1
\RequirePackage{ifpdf}
\ifpdf 
\documentclass[pdftex]{sigma}
\else
\documentclass{sigma}
\fi

\numberwithin{equation}{section}

\newtheorem{Theorem}{Theorem}[section]
\newtheorem*{Theorem*}{Theorem}

 { \theoremstyle{definition}

\newtheorem{Remark}[Theorem]{Remark} }

\def\phase{{\rm ph}}
\def\ZZ{{\mathbb Z}}
\def\NN{{\mathbb N}}
\def\eps{{\varepsilon}}
\def\wh{\widehat}
\def\bigO{{\cal O}}

\begin{document}

\allowdisplaybreaks

\renewcommand{\thefootnote}{}

\newcommand{\arXivNumber}{2402.14537}

\renewcommand{\PaperNumber}{075}
	
\FirstPageHeading

\ShortArticleName{McMahon-Type Asymptotic Expansions of the Zeros of the Coulomb Wave Functions}

\ArticleName{McMahon-Type Asymptotic Expansions\\ of the Zeros of the Coulomb Wave Functions\footnote{This paper is a~contribution to the Special Issue on Asymptotics and Applications of Special Functions in Memory of Richard Paris. The~full collection is available at \href{https://www.emis.de/journals/SIGMA/Paris.html}{https://www.emis.de/journals/SIGMA/Paris.html}}}

\Author{Amparo GIL~$^{\rm a}$, Javier SEGURA~$^{\rm b}$ and Nico M. TEMME~$^{\rm c}$}

\AuthorNameForHeading{A.~Gil, J.~Segura and N.M.~Temme}

\Address{$^{\rm a)}$~Departamento de Matem\'atica Aplicada y CC, de la Computaci\'on, ETSI Caminos,\\
\hphantom{$^{\rm a)}$}~Universidad de Cantabria, 39005 Santander, Spain}
\EmailD{\href{mailto:amparo.gil@unican.es}{amparo.gil@unican.es}}
\URLaddressD{\url{http://personales.unican.es/gila/}}

\Address{$^{\rm b)}$~Departamento de Matem\'aticas, Estadistica y Computaci\'on, Universidad de Cantabria,\\
\hphantom{$^{\rm b)}$}~39005 Santander, Spain}
\EmailD{\href{mailto:javier.segura@unican.es}{javier.segura@unican.es}}
\URLaddressD{\url{http://personales.unican.es/segurajj/}}

\Address{$^{\rm c)}$~Valkenierstraat 25, 1825BD Alkmaar, The Netherlands}
\EmailD{\href{mailto:nic@temme.net}{nic@temme.net}}

\ArticleDates{Received February 23, 2024, in final form August 07, 2024; Published online August 10, 2024}

\Abstract{We derive asymptotic expansions of the large zeros of the Coulomb wave functions and for those of their derivatives. The new expansions have the same form as the McMahon expansions of the zeros of the Bessel functions and reduce to them when a parameter is equal to zero. Numerical tests are provided to demonstrate the accuracy of the expansions.
}


\Classification{33C47; 33C10; 33C15; 41A60; 65D20; 65H05}

\begin{flushright}
\begin{minipage}{53mm}
\it Dedicated to Richard B. Paris,\\ the man who loved asymptotics
\end{minipage}
\end{flushright}

\section{Introduction}\label{sec:intro}
The Coulomb wave functions $F_\lambda(\eta,\rho)$ and $G_\lambda(\eta,\rho)$ are two linearly independent solutions of the differential equation
\[
\frac{{\rm d}^2 w}{{\rm d} \rho^2}+\left(1-\frac{2\eta}{\rho}-\frac{\lambda(\lambda+1)}{\rho^2}\right)w=0.
\]
This equation can be transformed into Kummer's differential equation
and the relation with the Kummer functions is
\begin{gather}
F_\lambda(\eta,\rho)=A {}_1F_1(\lambda+1-{\rm i}\eta,2\lambda+2;2{\rm i}\rho),\nonumber\\
G_\lambda(\eta,\rho)={\rm i}F_\lambda(\eta,\rho)+{\rm i}B U(\lambda+1-{\rm i}\eta,2\lambda+2,2{\rm i}\rho),\nonumber\\
A=\frac{|\Gamma(\lambda+1+{\rm i}\eta)|{\rm e}^{-\pi\eta/2-{\rm i}\rho}(2\rho)^{\lambda+1}}{2\Gamma(2\lambda+2)},\qquad
B={\rm e}^{\pi\eta/2+\lambda\pi {\rm i}-{\rm i}\sigma_\lambda(\eta)-{\rm i}\rho}(2\rho)^{\lambda+1},\nonumber\\
\sigma_\lambda(\eta)=\phase\, \Gamma(\lambda+1+{\rm i}\eta)\qquad \text{(the Coulomb phase shift)}.\label{eq:intro02}
\end{gather}

Coulomb wave functions and their zeros find application in various fields. In particular, they play an important role in atomic and nuclear physics where they contribute to understanding phenomena like electron scattering and bound-state properties (see, for example, \cite{Luna:2019:LBS}).

In physics applications, $\lambda$ usually has the values $0,1,2,\dots$, but here we assume that $\lambda$ and $\eta$ are real. The Kummer function in $F_\lambda(\eta,\rho)$ is analytic for all complex values of the parameters, unless $2\lambda+2=-1,-2,-3,\dots$. But because of the reciprocal gamma function in the quantity~$A$, the function $F_\lambda(\eta,\rho)$ is analytic at these points. If $\lambda>-1$, then the function $F_\lambda(\eta,\rho)$ disappears at $\rho=0$, which follows from the algebraic factor $\rho^{\lambda+1}$. Because of this factor, if~${\lambda\notin\ZZ}$, then~$F_\lambda(\eta,\rho)$ becomes multivalued with a branch cut along the negative axis. If we are interested in negative zeros of $F_\lambda(\eta,\rho)$, we can use Kummer's relation for the Kummer function, which in this case is
\[
{}_1F_1(\lambda+1-{\rm i}\eta,2\lambda+2;2{\rm i}\rho)={\rm e}^{2{\rm i}\rho}{}_1F_1(\lambda+1+{\rm i}\eta,2\lambda+2;-2{\rm i}\rho),
\]
and we can use the method for positive zeros by changing the sign of $\eta$.

The function $G_\lambda(\eta,\rho)$ becomes unbounded when $\rho\to0$, unless when $\lambda=\eta=0$, in which case
\begin{equation}\label{eq:intro04}
F_0(0,\rho)=\sin\rho,\qquad G_0(0,\rho)=\cos\rho.
\end{equation}

The asymptotic expansions of the large zeros of the Coulomb functions to be given in this paper are related with those of the $z$-zeros of the Bessel functions $J_\nu(z)$ and $Y_\nu(z)$ and their derivatives derived by McMahon \cite{McMahon:1894:RBF}. To explain what happens for the Bessel function $J_\nu(z)$, we use the well-known representation
\begin{equation}\label{eq:intro05}
J_\nu(z)=\sqrt{\frac{2}{\pi z}}(\cos\theta(\nu,z)\phi_\nu(z)-\sin\theta(\nu,z)\psi_\nu(z)),
\end{equation}
where $\theta(\nu,z)=z-\frac12\nu\pi-\frac14\pi$. For large values of $z$, we have $\phi_\nu(z)=1+\bigO\big(1/z^2\big)$, $\psi_\nu(z)=\bigO(1/z)$, and the complete asymptotic expansions of these functions are derived in \cite[Section~7.4]{Watson:1944:TTB} for the Hankel functions. For details, we refer to \cite[Section~10.17\,(i)]{Olver:2010:BFS}.

The relation in \eqref{eq:intro05} with asymptotic expansions of $\phi_\nu(z)$ and $\psi_\nu(z)$ was the starting point for McMahon \cite{McMahon:1894:RBF} to derive the asymptotic expansion of the large zeros of the Bessel function~$J_\nu(z)$ and similarly for related functions. For details on these expansions we refer to \cite[Section~10.21\,(vi)]{Olver:2010:BFS}.

For the Coulomb functions, the asymptotic expansions of the large zeros, and the methods to derive these expansions, have much in common with the expansions of the large zeros of the Bessel functions. This is not surprising because of the following observations.
\begin{enumerate}\itemsep=0pt
\item[(1)]
Firstly, because the formulas for the Coulomb functions given in the next section and used for deriving asymptotic expansions of the large zeros have the same analytical form as the one given in \eqref{eq:intro05} for $J_\nu(z)$.

\item[(2)]
Secondly, for $\eta=0$ we have \smash{$F_\lambda(0,\rho)=\sqrt{\pi\rho/2} J_{\lambda+\frac12}(\rho)$}, and we will verify in a special case that our expansion of the zeros of $F_\lambda(0,\rho)$ indeed become the McMahon expansion for the zeros of this Bessel function.
\end{enumerate}

What does surprise us is that we cannot find much information in the literature about the $\rho$-zeros of the Coulomb wave functions. Milton Abramowitz \cite{Abramowitz:AEC:1948} paid some attention to the zeros of $F_\lambda(\eta,\rho)$ in his article on asymptotic expansions of the Coulomb functions. His starting point is the same equation as ours, but he went further using an iteration method requiring function evaluations; we give more details in Section~\ref{sec:zeros}. Ikebe \cite{Ikebe:1975:ZRC} considered the zeros of $F_\lambda(\eta,\rho)$ and its derivative by computing eigenvalues of matrices following from recurrence relations of the Coulomb functions. See also \cite{Miyazaki:2001:EAC} for error analysis of this approach. Ball \cite{Ball:2000:ACZ} used a similar method for Bessel functions, Coulomb wave functions, and other special functions. This method is not based on function evaluations nor on asymptotic expansions, it requires eigenvalue computations of matrices. The early zeros can be computed efficiently. In fact, our method is efficient for the large zeros, and we show in examples how it performs for the first zeros.

We derive the asymptotic expansions for the zeros of both $F_\lambda(\eta,\rho)$ and $G_\lambda(\eta,\rho)$, and of their derivatives. The approach for these four functions is similar and the results are as simple as McMahon's expansion for zeros of the Bessel functions. A minor complication is that we need to solve a nonlinear equation, which can be done by standard numerical methods, although its solution can be expressed in terms of the Lambert $W$-function.

\section{Summary of used formulas}\label{sec:forms}

We summarise a set of formulas from the literature, see \cite[Chapter~31]{Temme:2015:AMI} and \cite[Section~33.11]{Thompson:2010:CWF}.

The following functions are important to describe the large $\rho$ asymptotics:
\begin{gather}
P_\lambda(\eta,\rho)=\sin(\theta_\lambda(\eta,\rho))F_\lambda(\eta,\rho)+\cos(\theta_\lambda(\eta,\rho))G_\lambda(\eta,\rho),\nonumber\\
Q_\lambda(\eta,\rho)=\cos(\theta_\lambda(\eta,\rho))F_\lambda(\eta,\rho)-\sin(\theta_\lambda(\eta,\rho))G_\lambda(\eta,\rho),\label{eq:forms01}
\end{gather}
where
\begin{equation}\label{eq:forms02}
\theta_\lambda(\eta,\rho)=\rho-\eta\ln(2\rho)-\frac12\lambda\pi+\sigma_\lambda(\eta),\qquad
\theta_\lambda^\prime(\eta,\rho)=1-{\eta}/{\rho}.
\end{equation}
Here, $\sigma_\lambda(\eta)$ is the Coulomb phase shift given in \eqref{eq:intro02}, and the prime denotes the derivative with respect to $\rho$.

By inverting \eqref{eq:forms01}, we have
\begin{gather}
F_\lambda(\eta,\rho)=\sin(\theta_\lambda(\eta,\rho))P_\lambda(\eta,\rho)+\cos(\theta_\lambda(\eta,\rho))Q_\lambda(\eta,\rho),\nonumber\\
G_\lambda(\eta,\rho)=\cos(\theta_\lambda(\eta,\rho))P_\lambda(\eta,\rho)-\sin(\theta_\lambda(\eta,\rho))Q_\lambda(\eta,\rho),\nonumber\\
F_\lambda^\prime(\eta,\rho)=\cos(\theta_\lambda(\eta,\rho))R_\lambda(\eta,\rho)+\sin(\theta_\lambda(\eta,\rho))S_\lambda(\eta,\rho),\nonumber\\
G_\lambda^\prime(\eta,\rho)=-\sin(\theta_\lambda(\eta,\rho))R_\lambda(\eta,\rho)+\cos(\theta_\lambda(\eta,\rho))S_\lambda(\eta,\rho),\label{eq:forms03}
\end{gather}
where
\[
R_\lambda(\eta,\rho)=P_\lambda(\eta,\rho)\theta_\lambda^\prime(\eta,\rho)+Q_\lambda^\prime(\eta,\rho),\qquad
S_\lambda(\eta,\rho)=P_\lambda^\prime(\eta,\rho)-Q_\lambda(\eta,\rho)\theta_\lambda^\prime(\eta,\rho).
\]
The functions $P_\lambda(\eta,\rho)$ and $Q_\lambda(\eta,\rho)$ can be written in terms of the Kummer $U$-functions. For details, we refer to the cited references.

The following asymptotic expansions follow from those of the $U$-function. We have for large values of $\rho$ the expansions
\begin{gather}
P_\lambda(\eta,\rho)\sim\sum_{k=0}^\infty \frac{p_k}{(2\rho)^k},\qquad
Q_\lambda(\eta,\rho)\sim\sum_{k=0}^\infty \frac{q_k}{(2\rho)^k},\nonumber\\
R_\lambda(\eta,\rho)\sim\sum_{k=0}^\infty \frac{r_k}{(2\rho)^k},\qquad
S_\lambda(\eta,\rho)\sim\sum_{k=0}^\infty \frac{s_k}{(2\rho)^k}.\label{eq:forms05}
\end{gather}
The coefficients of these expansions follow from simple recurrence relations. Initial values are
\[
p_0=1,\qquad q_0=0, \qquad r_0=1,\qquad s_0=0,
\]
and for $k=0,1,2,3,\dots$, we have
\begin{gather*}
(k+1)p_{k+1}=u_kp_k+v_k q_k,\qquad
(k+1)q_{k+1}=-v_kp_k+u_k q_k,\qquad
u_k=\eta(2k+1),\\
 v_k =k+k^2-\lambda^2-\lambda-\eta^2,\qquad
r_{k+1}=p_{k+1}-2\eta p_{k}-2kq_{k},\\
 s_{k+1}=-q_{k+1}+2\eta q_{k}-2kp_{k}.
\end{gather*}

\section{McMahon-type expansions of the zeros}\label{sec:zeros}

To start obtaining the expansion of the zeros of $F_\lambda(\eta,\rho)$, we look at the first line of \eqref{eq:forms03}. From the asymptotic expansions in \eqref{eq:forms05}, we see that
\[
P_\lambda(\eta,\rho)=1+\bigO(1/\rho),\qquad Q_\lambda(\eta,\rho)=\bigO(1/\rho),\qquad \rho\to\infty.
\]
Hence, for a large zero of $F_\lambda(\eta,\rho)$ the sine function in the first line of \eqref{eq:forms03} should be of order~$\bigO(1/\rho)$. We write
\begin{equation}\label{eq:zeros02}
\theta_\lambda(\eta,\rho)=n\pi+\delta,\qquad n\in\NN,
\end{equation}
and we will see that $\delta=\bigO(1/\rho)$. Using this form of $\theta_\lambda(\eta,\rho)$, we obtain
\[
 \sin(\theta_\lambda(\eta,\rho))=\cos(n\pi)\sin(\delta),\qquad
 \cos(\theta_\lambda(\eta,\rho))=\cos(n\pi)\cos(\delta),
\]
and from the first line of \eqref{eq:forms03}, we obtain
\begin{equation}\label{eq:zeros04}
\sin(\delta)P_\lambda(\eta,\rho)+\cos(\delta)Q_\lambda(\eta,\rho)=0,
\end{equation}
if we assume $F_\lambda(\eta,\rho)=0$.

Next, we try to find $\rho$ from equation \eqref{eq:zeros02}. We have
\begin{equation}\label{eq:zeros05}
\rho-\eta\ln\rho=\eta\ln 2+\frac12\lambda\pi-\sigma_\lambda(\eta)+n\pi +\delta,
\end{equation}
where $\delta$ is the small quantity introduced in \eqref{eq:zeros02} that must be found together with the zero $\rho$ for given values of $n$, $\eta$ and $\lambda$.

Let $\rho_0$ be the solution of the equation
\begin{equation}\label{eq:zeros06}
\rho_0-\eta\ln\rho_0=\eta\ln 2+\frac12\lambda\pi-\sigma_\lambda(\eta)+n\pi.
\end{equation}
Because the solution should satisfy $\rho_0=\bigO(n)$ for large $n$, we need to solve this equation for~${\rho_0>\eta}$.
Comparing \eqref{eq:zeros05} and \eqref{eq:zeros06}, we conclude
$
(\rho-\eta\ln\rho)-(\rho_0-\eta\ln\rho_0)=\delta$.
Also, with a new quantity $\eps$,
\begin{equation}\label{eq:zeros08}
\rho=\rho_0+\eps \quad \Longrightarrow \quad\delta=\eps-\eta\ln\left(1+\frac{\eps}{\rho_0}\right).
\end{equation}

With this $\rho_0$, $\delta$ and $\eps$ we try to find a solution $\rho$ of equation \eqref{eq:zeros04}.
First, assuming that~${\eps=\bigO(1/\rho_0)}$ for large $\rho_0$, we introduce the expansion
\begin{equation}\label{eq:zeros09}
\eps\sim\sum_{k=1}^\infty \frac{\eps_k}{\rho_0^k},\qquad \rho_0\to\infty.
\end{equation}
Using this expansion, we can obtain an expansion of $\delta$ in inverse powers of $\rho_0$ as well, and substitute this, with $\rho=\rho_0+\delta$ in \eqref{eq:zeros04}. We use the expansions of the functions $P_\lambda(\eta,\rho)$ and~${Q_\lambda(\eta,\rho)}$ given in \eqref{eq:forms05}, with $\rho=\rho_0+\eps$, and collect equal powers of $\rho_0$ to find the coefficients~$\eps_k$.

In this way, we obtain with $v_0=-\lambda^2-\lambda-\eta^2$,
\begin{equation}\label{eq:zeros10}
\eps_1=\frac12v_0,\qquad
\eps_2=\frac14\eta(3v_0+1),\qquad
\eps_3=\frac{1}{24}\big(22\eta^2v_0+17\eta^2-7v_0^2-6v_0\big).
\end{equation}
With these coefficients we find the wanted asymptotic expansion of the solution $\rho$, denoted by~$\rho_n$, of equation \eqref{eq:zeros05}
\begin{equation}\label{eq:zeros11}
\rho_n\sim \rho_0+
\frac{\eps_1}{\rho_0}+\frac{\eps_2}{\rho_0^2}+\frac{\eps_3}{\rho_0^3}+\cdots,\qquad n\to\infty.
\end{equation}

From numerical tests, we conclude that this gives indeed the approximation for the $n$-th positive zero.
Using the expansion in \eqref{eq:zeros11}, with $n=1$, $\lambda=2$, $\eta=\frac32$ we find ${\rho_0\doteq 9.186}$ and ${F_\lambda(\eta,\rho_{1})\doteq -0.0269}$. For $n=10$ and this expansion with the same $\lambda$, $\eta$ values, we find~${\rho_0\doteq39.65}$ and $F_\lambda(\eta,\rho_{10})\doteq 0.0000530$.
More tests are given in Section~\ref{sec:tests}.

\begin{Remark}\label{rem:rem01}
Abramowitz wrote in his paper \cite{Abramowitz:AEC:1948} on the asymptotics of the Coulomb wave functions our formula \eqref{eq:zeros04} in the form $\delta=-\arctan (Q_\lambda(\eta,\rho)/P_\lambda(\eta,\rho)$ and used \eqref{eq:zeros02} and \eqref{eq:forms02} to define the iteration
\[
\rho_{n,s}=\eta\ln(2\rho_{n,s-1})+\frac12\lambda\pi-\sigma_\lambda(\eta)+n\pi-
\arctan\frac{Q_\lambda(\eta,\rho_{n,s-1})}{P_\lambda(\eta,\rho_{n,s-1})},
\]
for $s=1,2,3,\dots$, where a starting value $\rho_{n,0}$ is needed, and the evaluation of the functions in the arctan-function.

For the analogues of our functions $P_\lambda(\eta,\rho)$ and $Q_\lambda(\eta,\rho)$, Abramowitz derived asymptotic expansions for large $\rho$, with $\lambda=0$, which are similar to those in the first line of \eqref{eq:forms05}. He did not use expansions of the Kummer $U$-function, but he derived the expansions using integral representations of functions related to the Coulomb functions, just as Hankel functions can be used for the asymptotics of the Bessel functions.

Abramowitz computed the first three zeros of $F_0(\eta,\rho)$, for a few $\eta$-values. See Table~\ref{table5} in Section~\ref{sec:tests} for a selection of these values. He did not give starting values but one may try $\rho_{n,0}=\rho_0$ defined in~\eqref{eq:zeros06}.
\end{Remark}

\begin{Remark}\label{rem:rem02}
When $\eta=0$, the $F$-Coulomb functions become $J$-Bessel functions. We have \smash{$F_\lambda(0,\rho)=\sqrt{\pi \rho/2}J_{\lambda+\frac12}(\rho)$} and \eqref{eq:zeros06} gives $\rho_0=\left(\frac12\lambda+n\right)\pi$. It is not difficult to verify that the first coefficients given in \eqref{eq:zeros10} become the first coefficients in McMahon's expansion of the zeros of the $J$-Bessel function, see \cite[Section~10.21\,(vi)]{Olver:2010:BFS}.
\end{Remark}

For the function $G_\lambda(\eta,\rho)$, we see in the second line of~\eqref{eq:forms01} that for a first approximation we have to use the zeros of the cosine function, and we change \eqref {eq:zeros02} and \eqref{eq:zeros06} by replacing~$n$ with~${n-\frac12}$.
This defines $\delta$ for this case and with this new $\theta_\lambda(\eta,\rho)$ the second line of \eqref{eq:forms03} becomes the same as in
 \eqref{eq:zeros04}. It follows that the coefficients $\eps_k$ in the expansion of $\eps$ in \eqref{eq:zeros09} are the same as those for the zeros of $F_\lambda(\eta,\rho)$, with different $\rho_0$ in the expansion. This property corresponds with McMahon's expansion for the zeros of the $Y$-Bessel function: the asymptotic expansions of the zeros of $J_\nu(x)$ and $Y_\nu(x)$ have the same coefficients, but the series have a~different large parameter.

 For the derivative $F_\lambda^\prime(\eta,\rho)$, we change \eqref {eq:zeros02} and \eqref{eq:zeros06} similarly, replacing $n$ with $n-\frac12$, and the equation for the zeros corresponding with \eqref{eq:zeros04} becomes
 \begin{equation}\label{eq:zeros13}
-\sin(\delta)R_\lambda(\eta,\rho)+\cos(\delta)S_\lambda(\eta,\rho)=0.
\end{equation}
This gives new coefficients for the expansion of $\eps$, and we write \eqref{eq:zeros08} as
\[
\wh\rho=\rho_0+\wh\eps \quad \Longrightarrow \quad\delta=\wh\eps-\eta\ln\left(1+\frac{\wh\eps}{\rho_0}\right),
\]
and the expansions corresponding with \eqref{eq:zeros09} and \eqref{eq:zeros11} in the form
\begin{equation}\label{eq:zeros15}
\wh\eps\sim\sum_{k=1}^\infty \frac{\wh\eps_k}{\rho_0^k},\qquad
\wh\rho_n\sim \rho_0+
\frac{\wh\eps_1}{\rho_0}+\frac{\wh\eps_2}{\rho_0^2}+\frac{\wh\eps_3}{\rho_0^3}+\cdots,\qquad n\to\infty.
\end{equation}
The first coefficients are, again with $v_0=-\lambda^2-\lambda-\eta^2$,
 \begin{equation}\label{eq:zeros16}
\wh\eps_1=\frac12v_0,\qquad
\wh\eps_2=\frac14\eta(3v_0-1),\qquad
\wh\eps_3=\frac{1}{24}\big(22\eta^2v_0-19\eta^2-7v_0^2+6v_0\big).
\end{equation}

Finally, we consider $G_\lambda^\prime(\eta,\rho)$. We take, as for $F_\lambda(\eta,\rho)$,
\eqref{eq:zeros02} and \eqref{eq:zeros06}. Hence, the expansions of the zeros of $G_\lambda^\prime(\eta,\rho)$ follow from \eqref{eq:zeros13} and have the same coefficients as in \eqref{eq:zeros16}.

We summarise for all four cases the choices of $\theta_\lambda(\eta,\rho)$, $\rho_0$ and the equation needed for obtaining the expansions.
\begin{enumerate}\itemsep=0pt
\item[(1)] $F_\lambda(\eta,\rho)$
\begin{gather*}
\theta_\lambda(\eta,\rho)=n\pi+\delta,\qquad n\in\NN,\\
\rho_0-\eta\ln\rho_0=\eta\ln 2+\frac12\lambda\pi-\sigma_\lambda(\eta)+n\pi,\qquad
\sin(\delta)P_\lambda(\eta,\rho)+\cos(\delta)Q_\lambda(\eta,\rho)=0.
\end{gather*}

\item[(2)] $G_\lambda(\eta,\rho)$
\begin{gather*}
\theta_\lambda(\eta,\rho) = \left(n-\frac12\right)\pi+\delta,\qquad n\in\NN,\\
\rho_0-\eta\ln\rho_0 = \eta\ln 2+\frac12\lambda\pi-\sigma_\lambda(\eta)+\left(n-\frac12\right)\pi,\\
\sin(\delta)P_\lambda(\eta,\rho) + \cos(\delta)Q_\lambda(\eta,\rho)=0.
\end{gather*}

\item[(3)] $F_\lambda^\prime(\eta,\rho)$
\begin{gather*}
\theta_\lambda(\eta,\rho) = \left(n-\frac12\right)\pi+\delta,\qquad n\in\NN,\\
\rho_0-\eta\ln\rho_0 = \eta\ln 2+\frac12\lambda\pi-\sigma_\lambda(\eta)+\left(n-\frac12\right)\pi,\\
-\sin(\delta)R_\lambda(\eta,\rho) + \cos(\delta)S_\lambda(\eta,\rho)=0.
\end{gather*}

\item[(4)] $G_\lambda^\prime(\eta,\rho)$
\begin{gather*}
\theta_\lambda(\eta,\rho) = n\pi+\delta,\qquad n\in\NN,\\
\rho_0-\eta\ln\rho_0 = \eta\ln 2+\frac12\lambda\pi-\sigma_\lambda(\eta)+n\pi,\qquad
-\sin(\delta)R_\lambda(\eta,\rho) + \cos(\delta)S_\lambda(\eta,\rho)=0.
\end{gather*}
\end{enumerate}

For all four cases, we have verified by using numerical calculations and graphs that the $n$ used in these relations
approximate the $n$th positive zero of the
Coulomb wave function or its derivative. Some of the tests performed are shown in the next section.
For $\lambda=\eta=0$, the first positive zeros have the proper $n$-value, which trivially follows from \eqref{eq:intro04}.

\section{Numerical tests}\label{sec:tests}

For testing the approximations obtained with the McMahon-type expansions $\big(\rho_n^{Mc}\big)$ for the zeros of Coulomb functions, we use the numerical method described in \cite{Segura:2010:zerSF} implemented in \textsc{Maple} with a large number of digits. Tables \ref{table1}, \ref{table2}, \ref{table3} and \ref{table4}
 show tests for the first $10$ zeros of the functions $F_{1.3}(2.1,\rho)$, $G_{1.3}(2.1,\rho)$, $F'_{1.3}(2.1,\rho)$ and $G'_{1.3}(2.1,\rho)$, respectively.
We use $6$ terms in the expansions \eqref{eq:zeros11} and \eqref{eq:zeros15}. High-accuracy numerical values of the zeros of the functions
obtained with the method described in \cite{Segura:2010:zerSF}, are given in the second column of the tables.
The third column in all tables shows the relative errors obtained in the comparisons. As can be seen, an accuracy close to $10^{-8}$ can be obtained for the largest considered zeros. As expected, the accuracy of the approximations improves as $n$ increases.

\begin{table}[th]\renewcommand{\arraystretch}{1.19}
$$
\begin{array}{rllll}
\hline
 n & \rho_n^{Mc} & \rho_n , \mbox{using the method in \cite{Segura:2010:zerSF}} & \mbox{Rel.\ error} \\
 \hline
 1 & \underline{9.2}8\ldots & 9.276226087098264 & 6.8 \times 10^{-4} \\
 2 & \underline{13.32}1\ldots & 13.32061436693835 & 5.3 \times 10^{-5} \\
 3 & \underline{17.049}4\ldots & 17.04925305758087 & 9.0 \times 10^{-6} \\
 4 & \underline{20.633}2\ldots & 20.63316305105047 & 2.3 \times 10^{-6} \\
 5 & \underline{24.1319}8\ldots & 24.13196399208639 & 7.5 \times 10^{-7} \\
 6 & \underline{27.5741}5\ldots & 27.57414717920683 & 2.9 \times 10^{-7} \\
 7 & \underline{30.97572}9\ldots & 30.97572598757761 & 1.3 \times 10^{-7} \\
 8 & \underline{34.34666}2\ldots & 34.34666006955555 & 6.1 \times 10^{-8} \\
 9 & \underline{37.69359}3\ldots & 37.69359261174668 & 3.2 \times 10^{-8} \\
 10 & \underline{41.02118}9\ldots & 41.02118854245900 & 1.7 \times 10^{-8} \\
\end{array}
$$

\vspace{-2mm}

\caption{Test for the McMahon-type approximations $\rho_n^{Mc}$ to the first $10$ zeros of the Coulomb function~$F_{1.3}(2.1,\rho)$.}\label{table1}
\end{table}

\begin{table}[th]\renewcommand{\arraystretch}{1.19}
$$
\begin{array}{rllll}
\hline
 n & \rho_n^{Mc} & \rho_n , \mbox{using the method in \cite{Segura:2010:zerSF}} & \mbox{Rel.\ error} \\
 \hline
 1 & \underline{6.9}5\ldots & 6.925107084382577 & 4.9 \times 10^{-3} \\
 2 & \underline{11.3}6\ldots & 11.35971565567721 & 1.6 \times 10^{-4} \\
 3 & \underline{15.209}4\ldots & 15.20913702648054 & 2.0 \times 10^{-5} \\
 4 & \underline{18.854}5\ldots & 18.85445602183751 & 4.3 \times 10^{-6} \\
 5 & \underline{22.3910}3\ldots & 22.39100849194709 & 1.2 \times 10^{-6} \\
 6 & \underline{25.8589}5\ldots & 25.85894221100473 & 4.6 \times 10^{-7} \\
 7 & \underline{29.2792}9\ldots & 29.27928968958546 & 1.9 \times 10^{-7} \\
 8 & \underline{32.66455}3\ldots & 32.66455053595783 & 8.7 \times 10^{-8} \\
 9 & \underline{36.022}8\ldots & 36.02279903910762 & 4.3 \times 10^{-8} \\
 10 & \underline{39.35957}2\ldots & 39.35957112638164 & 2.3 \times 10^{-8} \\
\end{array}
$$

\vspace{-2mm}

\caption{Test for the McMahon-type approximations $\rho_n^{Mc}$ to the first $10$ zeros of the Coulomb function~$G_{1.3}(2.1,\rho)$.}\label{table2}

\end{table}

\begin{table}[th!]\renewcommand{\arraystretch}{1.19}
$$
\begin{array}{rllll}
\hline
 n & \rho_n^{Mc} & \rho_n , \mbox{using the method in \cite{Segura:2010:zerSF}} & \mbox{Rel.\ error} \\
 \hline
 1 & \underline{6.}8\ldots & 6.740012285516214 & 2.0 \times 10^{-2} \\
 2 & \underline{11.3}4\ldots & 11.33586159146655 & 4.5 \times 10^{-4} \\
 3 & \underline{15.}2\ldots & 15.19947063325694 & 5.3 \times 10^{-5} \\
 4 & \underline{18.849}3\ldots & 18.84912765706333 & 1.1 \times 10^{-5} \\
 5 & \underline{22.3876}7\ldots & 22.38760195810186 & 3.2 \times 10^{-6} \\
 6 & \underline{25.8565}9\ldots & 25.85656409550572 & 1.1 \times 10^{-6} \\
 7 & \underline{29.2775}4\ldots & 29.27752955366132 & 4.6 \times 10^{-7} \\
 8 & \underline{32.66319}9\ldots & 32.66319220425298 & 2.1 \times 10^{-7} \\
 9 & \underline{36.0217}2\ldots & 36.02171734983164 & 1.0 \times 10^{-7} \\
 10 & \underline{39.3586}9\ldots & 39.35868838281058 & 5.5 \times 10^{-8} \\
\end{array}
$$

\vspace{-2mm}

\caption{Test for the McMahon-type approximations $\rho_n^{Mc}$ to the first $10$ zeros of the function $F'_{1.3}(2.1,\rho)$.}\label{table3}

\end{table}

\begin{table}[th]\renewcommand{\arraystretch}{1.15}
$$
\begin{array}{rllll}
\hline
 n & \rho_n^{Mc} & \rho_n , \mbox{using the method in \cite{Segura:2010:zerSF}} & \mbox{Rel.\ error} \\
 \hline
 1 & \underline{9.2}4\ldots & 9.226939712774167 & 2.0 \times 10^{-3} \\
 2 & \underline{13.30}8\ldots & 13.30627800305222 & 1.4 \times 10^{-4} \\
 3 & \underline{17.042}6\ldots & 17.04225058479286 & 2.3 \times 10^{-5} \\
 4 & \underline{20.62}9\ldots & 20.62896049608348 & 5.7 \times 10^{-6} \\
 5 & \underline{24.1291}8\ldots & 24.12914248690917 & 1.8 \times 10^{-6} \\
 6 & \underline{27.5721}3\ldots & 27.57211363372210 & 7.0 \times 10^{-7} \\
 7 & \underline{30.9741}9\ldots & 30.97418664616960 & 3.1 \times 10^{-7} \\
 8 & \underline{34.34545}7\ldots & 34.34545207910902 & 1.4 \times 10^{-7} \\
 9 & \underline{37.6926}2\ldots & 37.69261810059473 & 7.5 \times 10^{-8} \\
 10& \underline{41.02038}5\ldots & 41.02038500317911 & 4.1 \times 10^{-8} \\
\end{array}
$$

\vspace{-2mm}

\caption{Test for the McMahon-type approximations $\rho_n^{Mc}$ to the first $10$ zeros of the function $G'_{1.3}(2.1,\rho)$.
}\label{table4}

\end{table}
\begin{table}[th]\renewcommand{\arraystretch}{1.15}
$$
\begin{array}{lllll}
\hline
\qquad \eta & \qquad n=2 & \qquad n=3 \\
 \hline
 & (1) \underline{10.97}4 & (1) \underline{14.56}7 \\
 \eta=1.5 & (2) \underline{10.9733}6 & (2) \underline{14.56633}7 \\
 & (3){\bf 10.97335} & (3) {\bf 14.566335 } \\
\hline
 & (1) \underline{12.40}3 & (1) \underline{16.110} \\
 \eta=2 & (2) \underline{12.405}3 & (2)\underline{16.1104}7 \\
 & (3) {\bf 12.4052} & (3) {\bf 16.11044} \\
\hline
 & (1) \underline{13.78}6 & (1) \underline{17.59}6 \\
 \eta=2.5 & (2) \underline{13.78}85 & (2) \underline{17.595}4 \\
 & (3) {\bf 13.7879} & (3) {\bf 17.5953} \\

\hline
 & (1) \underline{15.13}0 & (1) \underline{19.03}3 \\
 \eta=3 & (2) \underline{15.13}49 & (2) \underline{19.035}6 \\
 & (3) {\bf 15.1335} & (3) {\bf 19.0352} \\

\end{array}
$$

\vspace{-2mm}

\caption{Comparison of values for the second and third zeros of $F_0(\eta,\rho)$ for a few $\eta$-values. (1) Values given in \cite{Abramowitz:AEC:1948}; (2) McMahon-type approximations; (3) Values obtained with the numerical method given in~\cite{Segura:2010:zerSF}.\looseness=-1}\label{table5}

\end{table}

\begin{figure}[th!]
\centering
\includegraphics[width=12.0cm]{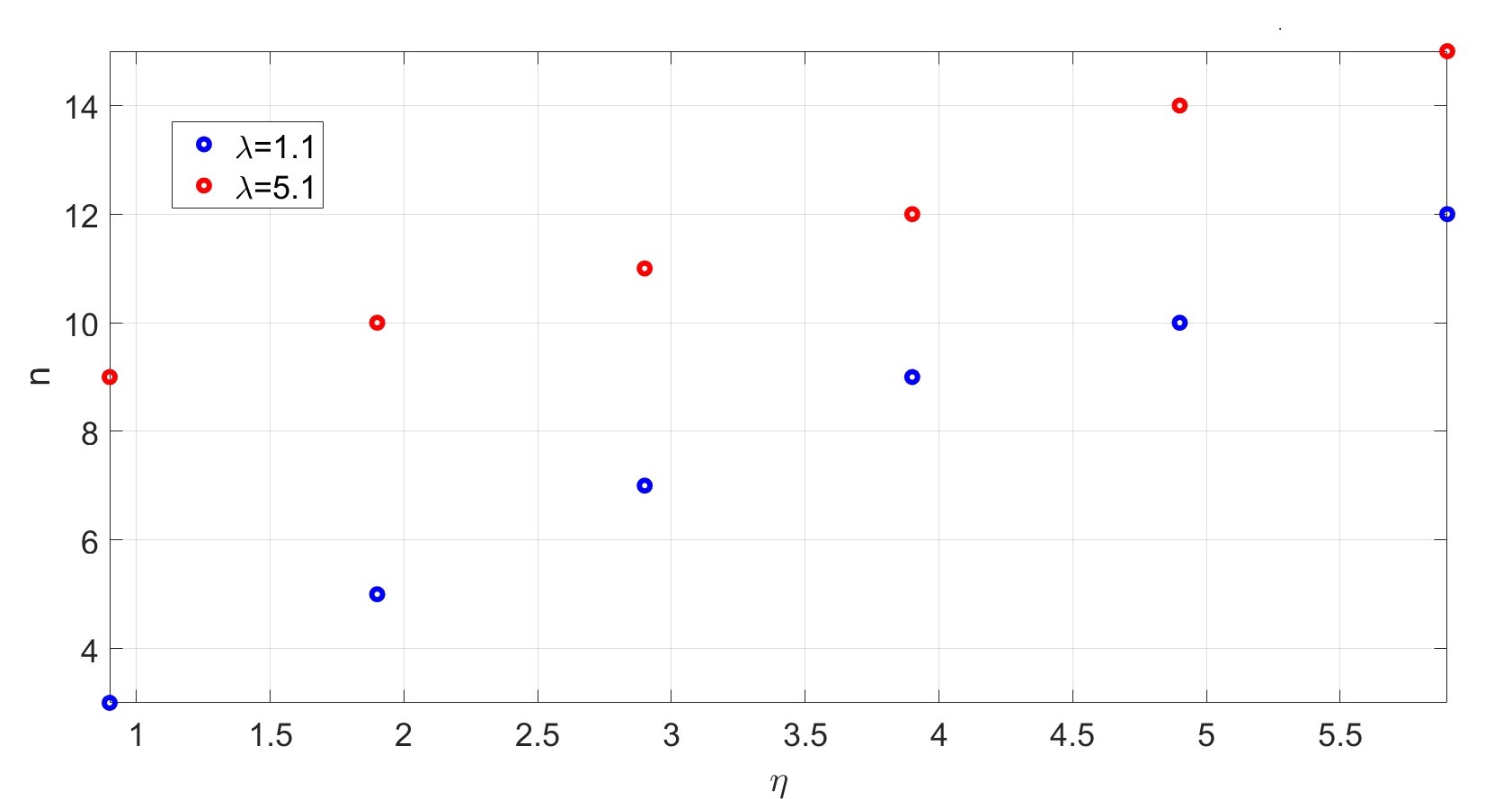}
\caption{Minimum value of $n$ for which the relative accuracy obtained with the McMahon-type
approximations $\rho_n^{Mc}$ to the zeros of $F_{\lambda}(\eta,\rho)$ is better than $10^{-6}$.}\label{fig:fig1}
\end{figure}

Additionally, a test of the influence of the parameters $\lambda$, $\eta$ on the accuracy of the approximations $\rho_n^{Mc}$
to the zeros of $F_{\lambda}(\eta,\rho)$, is shown in Figure \ref{fig:fig1}.
In this figure, we show the minimum value of $n$ for which the relative accuracy obtained with the McMahon-type
approximations is better than $10^{-6}$. The calculations have been made by fixing the value of $\lambda$ and varying the
values of $\eta$. The accuracy is checked, as before, using the numerical method given in \cite{Segura:2010:zerSF}. Results obtained for two different
values of $\lambda$ are shown for comparison. Similar results are obtained for the approximations to the zeros of the other functions.

As a final test, in Table \ref{table5}, we compare the McMahon-type
approximations with some of the values for the zeros of $F_0(\eta,\rho)$ appearing in the table given in Abramowitz's paper
\cite{Abramowitz:AEC:1948} for a~few~$\eta$-values. The values obtained with the numerical method \cite{Segura:2010:zerSF} are also shown.
 As can be seen, few discrepancies in the last digit for some of the values given in \cite{Abramowitz:AEC:1948} are found.

\subsection*{Acknowledgements}
We thank the referees for helpful and constructive remarks on an earlier version of the paper.
We acknowledge financial support from Ministerio de Ciencia e Innovaci\'on, Spain,
project PID2021-127252NB-I00 (MCIN/AEI/10.13039/501100011033/ FEDER, UE).

\pdfbookmark[1]{References}{ref}
\LastPageEnding
		

\begin{thebibliography}{99}
\footnotesize\itemsep=0pt

\bibitem{Abramowitz:AEC:1948}
Abramowitz M., Asymptotic expansions of {C}oulomb wave functions,
 \href{https://doi.org/10.1090/qam/28479}{\textit{Quart. Appl. Math.}} \textbf{7} (1949), 75--84.

\bibitem{Ball:2000:ACZ}
Ball J.S., Automatic computation of zeros of {B}essel functions and other
 special functions, \href{https://doi.org/10.1137/S1064827598339074}{\textit{SIAM~J. Sci. Comput.}} \textbf{21} (1999),
 1458--1464.

\bibitem{Ikebe:1975:ZRC}
Ikebe Y., The zeros of regular {C}oulomb wave functions and of their
 derivatives, \href{https://doi.org/10.2307/2005300}{\textit{Math. Comp.}} \textbf{29} (1975), 878--887.

\bibitem{Luna:2019:LBS}
Luna B.K., Papenbrock T., Low-energy bound states, resonances, and scattering
 of light ions, \href{https://doi.org/10.1103/PhysRevC.100.054307}{\textit{Phys. Rev.~C}} \textbf{100} (2019), 054307, 17~pages,
 \href{https://arxiv.org/abs/1907.11345}{arXiv:1907.11345}.

\bibitem{McMahon:1894:RBF}
Mcmahon J., On the roots of the {B}essel and certain related functions,
 \href{https://doi.org/10.2307/1967501}{\textit{Ann. of Math.}} \textbf{9} (1894), 23--30.

\bibitem{Miyazaki:2001:EAC}
Miyazaki Y., Kikuchi Y., Cai D., Ikebe Y., Error analysis for the computation
 of zeros of regular {C}oulomb wave function and its first derivative,
 \href{https://doi.org/10.1090/S0025-5718-00-01241-2}{\textit{Math. Comp.}} \textbf{70} (2001), 1195--1204.

\bibitem{Olver:2010:BFS}
Olver F.W.J., Maximon L.C., {B}essel functions, in N{IST} {H}andbook of
 {M}athematical {F}unctions, Cambridge University Press, Cambridge, 2010,
 215--286.

\bibitem{Segura:2010:zerSF}
Segura J., Reliable computation of the zeros of solutions of second order
 linear {ODE}s using a fourth order method, \href{https://doi.org/10.1137/090747762}{\textit{SIAM~J. Numer. Anal.}}
 \textbf{48} (2010), 452--469.

\bibitem{Temme:2015:AMI}
Temme N.M., Asymptotic methods for integrals, \textit{Ser. Anal.}, Vol.~6,
 \href{https://doi.org/10.1142/9195}{World Scientific Publishing}, Hackensack, NJ, 2014.

\bibitem{Thompson:2010:CWF}
Thompson I.J., {C}oulomb wave functions, in N{IST} {H}andbook of {M}athematical
 {F}unctions, Cambridge University Press, Cambridge, 2010, 741--756.

\bibitem{Watson:1944:TTB}
Watson G.N., A treatise on the theory of {B}essel functions, 2nd ed., \textit{Cambridge Math. Lib.}, Cambridge University Press, Cambridge, 1944.

\end{thebibliography}
\end{document}